\documentclass{amsart}

\usepackage{amsthm}
\usepackage{amsfonts}
\usepackage{amsmath}
\usepackage{amssymb}
\usepackage{mathrsfs}
\usepackage{fullpage}
\usepackage{stmaryrd}
\usepackage{hyperref}
\usepackage{verbatim}

\theoremstyle{plain}

\theoremstyle{definition}

\theoremstyle{remark}

\newcommand{\Begin}[2]{\begin{#1}\label{#2}}

\newcommand{\bPi}{\mathbf{\Pi}}
\newcommand{\bSigma}{\mathbf{\Sigma}}
\newcommand{\bDelta}{\mathbf{\Delta}}

\newcommand{\bbP}{\mathbb{P}}
\newcommand{\bbQ}{\mathbb{Q}}

\newcommand{\bbS}{\mathbb{S}}

\newcommand{\CA}{\mathcal{A}}
\newcommand{\CB}{\mathcal{B}}

\newcommand{\SCRL}{\mathscr{L}}

\newcommand{\proves}{\vdash}
\newcommand{\forces}{\Vdash}
\newcommand{\analytic}{{\bSigma_1^1}}
\newcommand{\lanalytic}{{\Sigma_1^1}}
\newcommand{\coanalytic}{{\bPi_1^1}}

\newcommand{\borel}{{\bDelta_1^1}}
\newcommand{\lborel}{{\Delta_1^1}}
\newcommand{\cantorspace}{{{}^\omega 2}}

\newcommand{\finNaturalSequence}{{{}^{<\omega}\omega}}
\newcommand{\KP}{\mathsf{KP}}

\newcommand{\F}{{F_{\omega_1}}}
\newcommand{\SR}{\mathrm{SR}}
\usepackage{enumerate}

\begin{document}

\title{Scott Ranks of Classifications of the Admissibility Equivalence Relation}

\author{William Chan}
\address{Department of Mathematics, University of North Texas, Denton, TX 76203}
\email{William.Chan@unt.edu}
\author{Matthew Harrison-Trainor}
\address{Department of Pure Mathematics, University of Waterloo, ON, Canada N2L 3G1}
\email{maharris@uwaterloo.ca}
\author{Andrew Marks}
\address{Department of Mathematics, UCLA, Los Angeles, CA 90095}
\email{marks@math.ucla.edu}

\begin{abstract}
Let $\SCRL$ be a recursive language. Let $S(\SCRL)$ be the set of $\SCRL$-structures with domain $\omega$. Let $\Phi : \cantorspace \rightarrow S(\SCRL)$ be a $\lborel$ function with the property that for all $x,y \in \cantorspace$, $\omega_1^x = \omega_1^y$ if and only if $\Phi(x) \approx_{\SCRL} \Phi(y)$. Then there is some $x \in \cantorspace$ so that $\SR(\Phi(x)) = \omega_1^x + 1$. 
\end{abstract}

\thanks{November 21, 2017. The first author was supported by NSF grant DMS-1703708. The second author was supported by a Banting fellowship. The third author was supported by NSF grant DMS-1500974.}

\maketitle


\section{Introduction}\label{Introduction}

The main equivalence relation of interest here is the countable admissible ordinal equivalence relation, denoted $\F$. It is an equivalence relation defined on $\cantorspace$ by $x \ \F \ y$ if and only if $\omega_1^x = \omega_1^y$. Recall that if $z \in \cantorspace$, then $\omega_1^z$ is the supremum of the collection of ordinals which are isomorphic to $z$-recursive well-orderings on $\omega$. Equivalently, $\omega_1^z$ is also the smallest $z$-admissible ordinal, i.e. the smallest ordinal height of a transitive model of $\KP$ containing $z$. The latter will be the more useful characterization here.

The equivalence relation $\F$ is important and can be meaningfully studied due to its connection with admissibility. A theorem of Sacks \cite{Countable-Admissible-Ordinals-and-Hyperdegrees} states that for any countable admissible ordinal $\alpha$, there is some $x \in \cantorspace$ so that $\omega_1^x = \alpha$. Therefore each equivalence class of $\F$ is associated with a countable admissible ordinal. $\F$ is a thin $\lanalytic$ equivalence relation with all equivalence classes $\borel$. (An equivalence relation $E$ is thin if and only if there is no perfect set of $E$-inequivalent elements.)

The topological Vaught's conjecture asserts that if the orbit equivalence relation of a Polish group acting continuously on a Polish space is thin, then it has countably many classes. Marker \cite{An-Analytic-Equivalence-Relation-Not-Arising} established a particular instance of this conjecture by showing that $\F$ is not an orbit equivalence relation of a continuous action of a Polish group on $\cantorspace$. This answered a question of Kechris. Becker \cite{Topological-Vaught-Conjecture-and-Mininal-Counterexamples} strengthened this by showing that $\F$ is not an orbit equivalence relation of a $\borel$ group action on $\cantorspace$. 

Suppose $E$ and $F$ are equivalence relations on Polish spaces $X$ and $Y$, respectively. $E$ is $\borel$ reducible to $F$, denoted $E \leq_\borel F$, if and only if there is a $\borel$ function $\Phi : X \rightarrow Y$ so that for all $a,b \in X$, $a \ E \ b \Leftrightarrow \Phi(a) \ F \ \Phi(b)$. $\borel$ reducibility is a common way of comparing the complexity of equivalence relations. 

Despite not being induced by a $\borel$ action of a Polish group, $\F$ is however $\borel$ reducible to an orbit equivalence relation of a continuous action of the Polish group $S_\infty$. Equivalence relations that are $\borel$ reducible to a continuous action of $S_\infty$ have a more model theoretic characterization:

Let $\SCRL$ be a countable language. Let $S(\SCRL)$ be the collection of $\SCRL$-structures with domain $\omega$. An equivalence relation $E$ on a Polish space $X$ is classifiable by countable structures if and only if there is a countable language $\SCRL$ so that $E \leq_\borel \approx_\SCRL$, where $\approx_\SCRL$ is the $\SCRL$-isomorphism relation defined on $S(\SCRL)$. A $\borel$ function $\Phi$ which witnesses this $\borel$ reducibility is called a classification of $E$ by $\SCRL$-structures. $\F$ is classifiable by countable structures by using the isomorphism relation of structures in the language of linear orderings. The following is an example of a classification function witnessing this.

Suppose $x \in \cantorspace$. An $x$-recursive pseudo-wellordering is an $x$-recursive linear ordering which is not a wellordering but has no $x$-hyperarithmetic descending sequences. Let $\eta$ be the order type of $\bbQ$. It is shown in \cite{Recursive-Pseudo-Well-Orderings} that $x$-recursive pseudo-wellorderings have ordertype $\omega_1^x (1 + \eta) + \rho$, for some $\rho < \omega_1^x$. (Also see \cite{Higher-Recursion-Theory}, Lemma III.2.2 (ii).) An $x$-Harrison linear ordering is an $x$-recursive linear ordering on $\omega$ of ordertype $\omega_1^x(1 + \eta)$. Note that the isomorphism type of an $x$-Harrison linear ordering is completely determined by $\omega_1^x$. The desired classification of $\F$ by linear orderings will be a map $\Phi$ sending $x$ to an $x$-Harrison linear ordering. Such a classification function exists which is $\borel$ follows from existence of a natural and uniform construction of Harrison linear orderings.

Although a recursive pseudo-wellordering can be constructed by a Barwise compactness argument, there is a more concrete natural construction found in \cite{On-the-Forms-Predicate-II} which is exposited in \cite{Higher-Recursion-Theory} Lemma III.2.1: The relation $R(x,y)$ if and only if $x$ is not $y$-hyperarithmetic is $\lanalytic$. As $\lanalytic$ sets are projections of recursive trees, one can use the relation $R$ to obtain a recursive tree $U$ in $2 \times \omega$ with the property that for all $y \in \cantorspace$, the tree $U^y$ is illfounded but has no $y$-hyperarithmetic path. It can be shown that the Kleene-Brouwer ordering of $U^y$ is a $y$-recursive pseudo-wellordering which can be modified to obtain a $y$-Harrison linear ordering. This procedure is $\lborel$ and gives the desired $\lborel$ classification $\Phi$. (See \cite{The-Countable-Admissible-Ordinal-Equivalence-Relation-pub} Section 3 for more details.)

Note that the Scott rank of $x$-Harrison linear orderings is $\omega_1^x + 1$. Therefore, in the example above, for all $x \in \cantorspace$, $\SR(\Phi(x)) = \omega_1^x + 1$. It is shown in \cite{The-Countable-Admissible-Ordinal-Equivalence-Relation-pub} Theorem 4.2 that any $\lborel$ classification $\Phi$ of $\F$ must send any real $x$ to an $x$-hyperarithmetic $\SCRL$-structure $\Phi(x)$ with $\SR(\Phi(x)) \geq \omega_1^x$. Recall that any $x$-hyperarithmetic $\SCRL$-structure $M$ on $\omega$ has Scott rank less than or equal to $\omega_1^x + 1$. If an $x$-hyperarithmetic structure $M$ has $\SR(M) \geq \omega_1^x$, then $M$ is said to have high Scott rank. 

\cite{The-Countable-Admissible-Ordinal-Equivalence-Relation-pub} Remark 4.3 raised the question of whether every $\lborel$ classification $\Phi$ has the property that for all $x$, $\SR(\Phi(x)) = \omega_1^x + 1$. An attempt to compel the Scott rank to take the highest possible value in \cite{The-Countable-Admissible-Ordinal-Equivalence-Relation-pub} Remark 4.3 fails due to the existence of $x$-recursive structures of Scott rank $\omega_1^x$. This suggests the following which is the main question of the paper:

\Begin{question}{SR classification high not highest rank}
Is there a recursive language $\SCRL$ and $\lborel$ function $\Phi : \cantorspace \rightarrow S(\SCRL)$ so that for all $x,y \in \cantorspace$, $x \ \F \ y \Leftrightarrow \Phi(x) \approx_\SCRL \Phi(y)$ with the property that for all $x \in \cantorspace$, $\SR(\Phi(x)) = \omega_1^x$?
\end{question}

Hyperarithmetic structures of Scott rank $\omega_1^\emptyset$ are more difficult to produce than computable structures of Scott rank $\omega_1^\emptyset + 1$. Makkai \cite{An-Example-Concerning-Scott-Heights} constructed the first such hyperarithmetic example. Knight and Millar \cite{Computable-Trees-of-Scott-Rank} produced recursive structures of Scott rank $\omega_1^\emptyset$. (Other examples can be found in \cite{Some-New-Computable-Structures-of-High-Rank}.) These constructions are much more intricate and have some nonuniform aspects due to the use of Barwise or Barwise-Kreisel compactness.

The main result of the paper is the following negative answer to Question \ref{SR classification high not highest rank}:
\\*
\\*
\noindent \textbf{Theorem \ref{classification has model high but not highest scott rank}.} Let $\SCRL$ be a recursive language and let $\Phi : \cantorspace \rightarrow S(\SCRL)$ be a $\lborel$ function so that $x \ \F \ y \Leftrightarrow \Phi(x) \approx_\SCRL \Phi(y)$. Then there is some $x \in \cantorspace$ so that $\SR(\Phi(x)) = \omega_1^x + 1$. 
\\*
\\*\indent (Becker has informed the authors that this theorem follows from a more general result in \cite{Strange-Structures-from-Computable-Model-Theory} which is however proved by different methods in $\mathsf{ZFC}$ augmented by $\bSigma_2^1$-determinacy. All results in this paper are implicitly proved from $\mathsf{ZFC}$.)

Theorem \ref{classification has model high but not highest scott rank} can be interpreted to say that there is no construction of a recursive structure of Scott rank $\omega_1^\emptyset$ which is natural enough in the sense that the construction can be relativized to any real $x$ to produce an $x$-recursive structure of Scott rank $\omega_1^x$ whose isomorphism type depends solely on the $x$-recursive ordinals. The notion of a $\lborel$ classification function for $\F$ is used to formalize this idea of naturality. This approach is similar to a questin of Martin concerning the unnaturalness of intermediate degrees: A well known result in recursion theory states that there is a degree strictly between $[\emptyset]_T$ and $[\emptyset']_T$. Martin asked whether there are any definable procedures that take an $X \in \mathcal{D}$, the set of Turing degrees, and return a Turing degree between $X$ and its jump. This is formalized by asking in $\mathsf{ZF + AD}$ whether there is any function $\Phi : \mathcal{D} \rightarrow \mathcal{D}$  with the property that for all $X \in \mathcal{D}$, $X <_T \Phi(X) <_T X'$. See \cite{A-Classification-of-Jump-Operators} and \cite{Martin-Conjecture-Arithmetic-Equivalence-Countable-Borel} for more on various forms of Martin's conjecture.

The proof of Theorem \ref{classification has model high but not highest scott rank} is similar to how \cite{On-the-Number-of-Countable-Models} Theorem 4.2 (or \cite{Bounds-on-Weak-Scattering} Corollay 6.2) shows that every counterexample to Vaught's conjecture $\tau$ has a model $M$ so that $\SR(M) = \omega_1^M + 1$. This is done by producing an illfounded end-extension of a $\Sigma_2$-admissible set which has an appropriate model of $\tau$ and has enough $\Sigma_1$-absoluteness to conclude that isolating formulas of the original $\Sigma_2$-admissible set are still isolating formulas in the illfounded end-extension.

An important feature of the proof of the above result for counterexamples to Vaught's conjecture is access to the complete theory and types of the desired model in small admissible sets even when the model does not exist in that admissible set. In the setting of this paper, access to the fragment, complete theory, and types of the desired model as well as the sufficient definability of these objects within the appropriate admissible set are obtained using the $\lborel$ classification function and the Solovay product forcing lemma for a suitable class forcing. This seems to be similar to ideas used in \cite{Computable-Structures-in-Generic-Extensions}.

Becker \cite{Strange-Structures-from-Computable-Model-Theory} has also considered the unnaturalness of $x$-recursive structures of Scott rank $\omega_1^x$. \cite{Strange-Structures-from-Computable-Model-Theory} shows that under $\mathsf{ZF + DC + AD}$, if $\SCRL$ is a countable language and $\mathcal{F}$ is a family of $\aleph_1$ many isomorphism types of $\SCRL$-structures, then there is a $z \in \cantorspace$ so that for all $x \geq_T z$, if an $x$-recursive $\SCRL$-structure $M$ has an isomorphism type in $\mathcal{F}$, then $\SR(M) \neq \omega_1^x$. Given any $\lborel$ classification $\Phi$ of $\F$, Theorem \ref{classification has model high but not highest scott rank} follows by applying Becker's result to the $\lanalytic$ family determined by the range of $\Phi$. \cite{Strange-Structures-from-Computable-Model-Theory} mentions that it is not known whether the $\analytic$ version of Becker's result holds in $\mathsf{ZFC}$ and that the methods of \cite{Strange-Structures-from-Computable-Model-Theory} to prove the $\analytic$ version require a determinacy assumption stronger than $\coanalytic$-determinacy but weaker than $\bPi_2^1$-determinacy.

\section{Basics}\label{Basics}

The results of the paper are proved in $\mathsf{ZFC}$. As customary in set theory, the real universe is denoted by $V$, which can be understood as some fixed model of $\mathsf{ZFC}$ where the results of the paper are being derived. Frequently concepts will be viewed from various different models of set theory. If $M$ is a model of set theory and $A$ is some notion given by a formula, $A^M$ will indicate the relativization of the definition of $A$ within the model $M$.

Let $\KP$ denote Kripke-Platek set theory with the infinity axiom, which can be formulated in any language $\mathscr{J}$ consisting of a distinguished binary relation symbol $\dot \in$ and possibly other symbols. $\KP$ is a weak axiom system for set theory. Its distinguishing axiom schemes are $\Delta_1$-separation and $\Sigma_1$-collection. Let $\Sigma_2$-$\KP$ be the axiom system extending $\KP$ by the axiom schemes of $\Delta_2$-separation and $\Sigma_2$-collection. 

An admissible set is a transive model of $\KP$. A $\Sigma_2$-admissible set is a transitive model of $\Sigma_2$-$\KP$. See \cite{Admissible-Sets-and-Structures} and \cite{Admissible-Sets} for more information about $\KP$ and admissibility.

Let $\mathrm{ON}$ denote the class of ordinals. If $\CA$ is some $\mathscr{J}$-structure satisfying $\KP$ where $\mathscr{J}$ is a language consisting of a distinguish binary relation symbol $\dot \in$, then $\mathrm{WF}(\CA)$ is the substructure of elements of $\CA$ which are $\dot\in^\CA$-well-founded in the real world $V$. Via the Mostowski collapse, one can always assume $(\mathrm{WF}(\CA), \dot\in^\CA)$ is a transitive set. The ordinal height of $\CA$ is $\mathrm{ON} \cap \mathrm{WF}(\CA) = \mathrm{ON} \cap \CA$, where it is assumed that $\mathrm{WF}(\CA)$ is a transitive set. 

For the rest of the paper, assume that $\omega$ belongs to the transitive closure of the well-founded part of any model of $\KP$.

\Begin{definition}{admissible ordinals}
Let $x \in \cantorspace$. $\alpha \in \mathrm{ON}$ is an $x$-admissible ordinal if and only if there is an admissible set $\CA$ so that $x \in \CA$ and $\alpha = \mathrm{ON} \cap \CA$. That is, $\alpha$ is the ordinal height of some admissible set containing $x$.  

The least $x$-admissible ordinal is denoted $\omega_1^x$. 
\end{definition}

\Begin{fact}{constructibility}
Let $\CA \models \KP$. There is a $\Delta_1$ function taking elements of $\alpha \in \mathrm{ON}^\CA$ to $L_\alpha$, the segment of G\"odel constructible hierarchy. $L^\CA$ is a $\Delta_1$ class in $\CA$ and $L^\CA \models \KP$. 

These results hold for the relativized G\"odel hierarchy.
\end{fact}

\Begin{definition}{hyp reals}
Let $x \in \cantorspace$. $y \in \cantorspace$ is an $x$-hyperarithmetic real if and only if $y$ belongs to every admissible set containing $x$. 
\end{definition}

A basic fact of descriptive set theory is that $y$ is $x$-hyperarithmetic if and only if $y$ is $\Delta_1^1(x)$. 

\Begin{fact}{smallest admissible}
Let $x \in \cantorspace$. $L_{\omega_1^x}[x]$ is the smallest $x$-admissible set under $\subseteq$. Hence the $x$-hyperarithmetic reals are exactly $\cantorspace \cap L_{\omega_1^x}[x]$. 
\end{fact}

\Begin{fact}{truncation lemma}
(Truncation lemma) If $\CA \models \KP$, then $\mathrm{WF}(\CA) \models \KP$. Therefore assuming that the well-founded part is transitive, $\mathrm{WF}(\CA)$ is an admissible set.
\end{fact}

\begin{proof}
See \cite{Admissible-Sets-and-Structures} Lemma II.8.4.
\end{proof}

The following result of Sacks gives an important characterization of countable admissible ordinals.

\Begin{fact}{sacks theorem}
(\cite{Countable-Admissible-Ordinals-and-Hyperdegrees}) If $\alpha$ is a countable admissible ordinal, then there is an $x \in \cantorspace$ so that $\omega_1^x = \alpha$. 
\end{fact}

This result can be proved using infinitary logic in countable admissible fragments as shown in \cite{Admissible-Sets}. These methods were used to study $\F$ in \cite{The-Countable-Admissible-Ordinal-Equivalence-Relation-pub} and will again be used in the arguments of this paper. The main tool for this approach is the Jensen's model existence theorem. 

\Begin{lemma}{jensen model existence theorem}
(Jensen's model existence theorem) Let $\CA$ be an admissible set. Let $\mathscr{J}$ be a language which is $\Delta_1$ definable in $\CA$ and contains a distinguished binary relation symbol $\dot \in$ and constant symbols $\hat{a}$ for each $a \in \CA$. Let $H$ be a consistent theory in the countable admissible fragment $(\mathscr{J}_{\infty\omega})^{\CA}$ associated with $\CA$ which is $\Sigma_1$ definable in $\CA$ and contains the following sentences:
\begin{enumerate}[(I)]
\item Extensionality.

\item For each $a \in \CA$, ``$(\forall v)(v \dot \in \hat{a} \Leftrightarrow \bigvee_{z \in a} v = \hat{z})$''.
\end{enumerate}
Then there is a $\mathscr{J}$-structure $\CB \models H$ so that $\mathrm{WF}(\CB)$ is transitive, $\CB$ end extends $\CA$, and $\mathrm{ON} \cap A = \mathrm{ON} \cap B$. 

If $\CA$ is a $\Sigma_2$-admissible set and the theory $H$ is $\Sigma_2$ definable in $\CA$, then the same conclusion holds. 
\end{lemma}

\begin{proof}
See \cite{Admissible-Sets} Section 4, Lemma 11 or \cite{Jensen-Model-Existence-Theorem}. Recall that $\CB$ is an end extension of $\CA$ if and only if $\CA \subseteq \CB$ and for all $x \in \CA$, $\{y \in \CA : y \dot\in^\CA x\} = \{y \in \CB : y \dot\in^\CB x\}$. 
\end{proof}

Fact \ref{sacks theorem} was originally proved by Sacks using a class forcing over countable admissible sets. For some properties concerning constructibility, the approach by forcing will be useful. A simple class forcing of Steel will be used. The following presents the definitions and basic properties. See \cite{Forcing-with-Tagged-Trees} for more details.

\Begin{definition}{steel forcing}
(Steel's forcing with tagged trees; see \cite{Forcing-with-Tagged-Trees}) Let $\CA$ be a countable model of $\KP$. Let $\infty$ be some symbol formally defined to be larger than all ordinals of $\CA$. Let $\bbS$ be the forcing consisting of $(T,h)$ where $T$ is a finite tree on $\omega$ and $h : T \rightarrow \mathrm{ON}^\CA \cup \{\infty\}$, with the property that for all $s,t \in T$ with $s \subseteq t$, $h(t) < h(s)$ or $h(s) = h(t) = \infty$. If $p,q \in \bbS$ and $p = (T_p,h_p)$ and $q = (T_q,h_q)$, then $p \leq_\bbP q$ if and only if $T_p \supseteq T_q$ and $h_p \supseteq h_q$. Let $1_\bbS = (\emptyset,\emptyset)$. The forcing relation $p \forces \varphi$, as a relation ranging over $p \in \bbS$ and ranked sentences $\varphi$ (see \cite{Forcing-with-Tagged-Trees}), is a $\Delta_1$ relation in $\CA$.

There are $\bbS$-names $\dot T, \dot h \in \CA$ so that for any $G \subseteq \bbS$ which is $\bbS$-generic over $\CA$, $\dot T[G] = \bigcup_{p \in G} T_p$ and $\dot h[G] = \bigcup_{p \in G} h_p$. Note that $\dot T[G]$ is a tree on $\omega$. When $G \subseteq \bbS$ is $\bbS$-generic over $\CA$, $\CA[\dot T[G]] \models \KP$. (However $\CA[G]$ is not a model of $\KP$.) Therefore, $\omega_1^{\dot T[G]} = \CA \cap \mathrm{ON}^V$.
\end{definition}

\Begin{definition}{infinitary logic definitions}
If $\SCRL$ is a language, then let $\SCRL_{\omega\omega}$ denote the set of first order $\SCRL$-formulas and $\SCRL_{\infty\omega}$ denote the class of infinitary formulas in the language $\SCRL$. In $\KP$, $\SCRL_{\infty\omega}$ is a $\Delta_1$ class. The satisfaction relation between $\SCRL$-structures and formulas of $\SCRL_{\infty\omega}$ is also $\Delta_1$ in $\KP$. A subset $\mathcal{F} \subseteq \SCRL_{\infty\omega}$ is an fragment if it has the closure properties of \cite{Admissible-Sets-and-Structures} Definition III.2.1. (See \cite{Admissible-Sets-and-Structures} Chapter III for more information about the syntax and semantics of $\SCRL_{\infty\omega}$.)
\end{definition}

\Begin{definition}{types definition}
Let $\SCRL$ be a recursive language, $\mathcal{F} \subseteq \SCRL_{\infty\omega}$ be a fragment, $T \subseteq \mathcal{F}$ be a theory, and $M$ be an $\SCRL$-structure. Then $S_n^{\mathcal{F}}(T)$ is the collection of all complete $n$-types of $T$ in the fragment $\mathcal{F}$.

For each $\varphi \in \mathcal{F}$ with $n$ many free variables, let $[\varphi]_{\mathcal{F}}^T = \{p \in S_n^\mathcal{F}(T) : \varphi \in p\}$. The topology on $S_n^\mathcal{F}(T)$ is generated by $[\varphi]_\mathcal{F}^T$ as basic open sets, where $\varphi$ ranges over formulas in $\mathcal{F}$ with $n$ free variables. A type $p \in S_n^\mathcal{F}(T)$ is an isolated type if $\{p\}$ is an open set. A type which is not isolated is sometimes called a nonprincipal type. A formula $\varphi$ is an isolating formula if and only $[\varphi]_\mathcal{F}^T$ is a singleton. That is, for all $\psi \in \mathcal{F}$ with $n$ free variables, $T \proves (\forall \bar{x})(\varphi \Rightarrow \psi)$ or $T \proves (\forall \bar{x})(\varphi \Rightarrow \neg\psi)$. 

If $\bar{a}$ is a tuple from $M$ of length $n$, then $\mathrm{tp}^M_{\mathcal{F}}(\bar{a})$ is the complete $n$-type consisting of the formulas of $\mathcal{F}$ satisfied by $\bar{a}$. 
\end{definition}

\Begin{definition}{theory of model and Scott rank}
The following definition and properties can be formalized and proved in $\KP$. 

Let $\SCRL$ be some language and let $M \in S(\SCRL)$. By $\Sigma_1$-recursion, the functions $\SCRL^M_\alpha$ and $T^M_\alpha$ are defined as follows:
\begin{enumerate}[$\bullet$]
\item Let $\SCRL^M_0 = \SCRL_{\omega\omega}$.

\item If $\alpha$ is a limit ordinals, then let $\SCRL^M_\alpha = \bigcup_{\beta < \alpha} \SCRL^M_\beta$. 

\item For any $\alpha$, let $T_\alpha^M$ be the complete theory of $M$ in the fragment $\SCRL^M_\alpha$. 

\item Let $\SCRL^M_{\alpha + 1}$ be the least fragment $\mathcal{F}$ extending $\SCRL^M_\alpha$ containing $\bigwedge p$ for each nonprincipal $p \in S_n^{\SCRL^M_\alpha}(T^M_\alpha)$ realized by some tuple in $M$. 
\end{enumerate}
The functions $\alpha \mapsto T^M_\alpha$ and $\alpha \mapsto \SCRL^M_\alpha$ are $\Sigma_1$-functions on the $\Delta_1$ class of ordinals. Hence these two functions are $\Delta_1$. Note that if $M, N \in S(\SCRL)$ and $M \approx_\SCRL N$, then $\SCRL^M_\alpha = \SCRL_\alpha^N$ and $T^M_\alpha = T^N_\alpha$ for all ordinals $\alpha$. 

Let $\SCRL^M_\infty = \bigcup_{\alpha \in \mathrm{ON}} \SCRL^M_\alpha$ and $T^M_\infty = \bigcup_{\alpha \in \mathrm{ON}} T_\alpha^M$. For any $\beta$, a formula of quantifer rank $\alpha < \beta$ belongs to $\SCRL^M_\beta$ or $T^M_\beta$ if it already belonged to $\SCRL^M_\alpha$ or $T^M_\alpha$. This can be used to show that $\SCRL^M_\infty$ and $T^M_\infty$ are $\Delta_1$.  

The Scott rank of $M$, denoted $\SR(M)$, is the smallest ordinal $\alpha$ so that $M$ is the atomic model of $T_\alpha^M$. 
\end{definition}

\Begin{fact}{bounds on Scott rank}
Let $\SCRL$ be a recursive language. Let $M \in S(\SCRL)$. $\SR(M) \leq \omega_1^M + 1$. 
\end{fact}

\begin{proof}
See \cite{Scott-Sentences-and-Admissible-Sets}.
\end{proof}

\section{Scott Ranks of Classifications}\label{Scott Ranks of Classifications}

Definition \ref{theory of model and Scott rank} shows that $T^M_\alpha$ and $\SCRL^M_\alpha$ can be defined for any structure $M$ within any model $\CA \models \KP$ containing $M$. Suppose $\sigma$ is a countable admissible ordinal and $x \in \cantorspace$ is such that $\omega_1^x = \sigma$. Let $\Phi$ be a $\Delta_1^1$ classification of $\F$ by $\SCRL$-structures. Generally, $\Phi(x) \notin L_\sigma$, where $L_\sigma$ is the $\sigma^\text{th}$-level of G\"odel constructible hierarchy and the smallest admissible set of height $\sigma$. Later, it will desirable to have $\SCRL^{\Phi(x)}_{\alpha}$ and $T^{\Phi(x)}_{\alpha}$ for $\alpha \leq \sigma$ either belong to $L_\sigma$ or are $\Delta_1$-definable in $L_\sigma$, even though the structure $\Phi(x)$ does not belong to $L_\sigma$.

 The Solovay product forcing lemma for class forcing will be very useful for showing that the relevant theory of a model that does not exist in $L_\sigma$ is actually definable in $L_\sigma$ under certain circumstances. The Solovay product lemma states that elements that belong to two mutually generic extensions actually already belong to the ground model.

\Begin{lemma}{solovay product lemma}
(Solovay product lemma) Let $\CA \models \KP$ . Let $\bbP$ be a $\Delta_1$-definable forcing in $\CA$. (This means that $\bbP$ and $\leq_\bbP$ are $\Delta_1$ definable.) Assume that $p \forces_{\bbP} \varphi$ is a $\Delta_1$ relation in arguments $p \in \bbP$ and ranked sentences $\varphi$. Let $\bbP \times \bbP$ denote the product forcing. Let $G,H \subseteq \bbP$ be $\bbP$-generic filters over $\CA$ such that $G \times H$ is a $\bbP \times \bbP$-generic filter over $\CA$. Then $\CA[G] \cap \CA[H] = \CA$.
\end{lemma}

\begin{proof}
Suppose $\CA[G] \cap \CA[H] \neq \emptyset$. Let $z \in (\CA[G] \cap \CA[H]) \setminus \CA$ be of minimal rank. Hence $z \subseteq \CA$. There are $\bbP$-names $\sigma$ and $\tau$ so that $z = \tau[G] = \sigma[H]$. By the forcing theorem, there is some $(p,q) \in \bbP \times \bbP$ so that $(p,q) \forces_{\bbP\times\bbP} \tau = \sigma$, where here $\tau$ and $\sigma$ are considered as $\bbP\times\bbP$-names which yield the result of the original $\bbP$-names $\tau$ and $\sigma$, respectively, evaluated using the left and right $\bbP$-generic filters, respectively, derived from $\bbP\times\bbP$-generic filters.

The claim is that for all $x \in \CA$, $p \forces_\bbP \check x \in \tau$ or $p \forces_\bbP \check x \notin \tau$: To see this, assume not. There is some $x \in \CA$ and some $p_0,p_1 \leq_\bbP p$ so that $p_0 \forces_{\bbP} \check x \in \tau$ and $p_1 \forces \check x \notin \tau$. Without loss of generality, suppose that $\CA[H] \models x \notin \sigma[H]$. Then find some $q' \leq_\bbP q$ so that $q' \forces_\bbP \check x \notin \sigma$. Let $G',H' \subseteq \bbP$ be $\bbP$-generic filters over $\CA$ so that $G'\times H'$ is a $\bbP\times\bbP$-generic over $\CA$ and $(p_0,q') \in G'\times H'$. By the forcing theorem, $x \in \tau[G']$ and $x \notin \sigma[H']$. Hence $\CA[G'\times H'] \models \tau[G'] \neq \sigma[H']$. But $(p_0,q') \leq_{\bbP\times\bbP} (p,q)$ and $(p,q) \forces_{\bbP\times\bbP} \tau = \sigma$. Contradiction. 

Let $\alpha$ be the rank of $z$. Let $\CA_\alpha$ denote the elements of $\CA$ of rank less than $\alpha$. Then $z = \{x \in \CA_\alpha : p \forces_{\bbP} \check x \in \tau\}$. $z \in \CA$ by $\Delta_1$-separation. This contradicts the earlier assumption that $z \notin \CA$. 
\end{proof}

\Begin{corollary}{exists in extension same omega1CK}
Let $\sigma$ be a countable admissible ordinal. $\bigcap_{\omega_1^x = \sigma} L_\sigma[x] = L_\sigma$.
\end{corollary}

\begin{proof}
Let $\bbS$ be the Steel's tagged tree forcing. Let $G,H$ be $\bbS$-generic filters over $L_\sigma$ so that $G \times H$ is $\bbS\times\bbS$-generic over $L_\sigma$. Let $a = \dot T[G]$ and $b = \dot T[H]$. $\omega_1^a = \omega_1^b = \sigma$. Then $\bigcap_{\omega_1^x = \sigma} L_\sigma[x] \subseteq L_\sigma[a] \cap L_\sigma[b] \subseteq L_\sigma[G] \cap L_\sigma[H] = L_\sigma$ using Lemma \ref{solovay product lemma}.
\end{proof}

In the following, let $\sigma$ be an admissible ordinal. Let $\SCRL$ be a recursive language. Let $\Phi : \cantorspace \rightarrow S(\SCRL)$ be a $\lborel$ classification of $\F$ by $\SCRL$-structures, i.e. $\omega_1^x = \omega_1^y$ if and only if $\Phi(x) \approx_\SCRL \Phi(y)$.

As mentioned above, since $\Phi(x) \approx_\SCRL \Phi(y)$ if and only if $\omega_1^x = \omega_1^y$, $\SCRL^{\Phi(x)}_\alpha = \SCRL^{\Phi(y)}_\alpha$ and $T^{\Phi(x)}_\alpha = T^{\Phi(y)}_\alpha$ for all ordinals $\alpha$ whenever $\omega_1^x = \omega_1^y$. Therefore, one may define $\SCRL^\sigma_\alpha$ and $T^\sigma_\alpha$ to be $\SCRL^{\Phi(x)}_\alpha$ and $T^{\Phi(x)}_\alpha$, respectively, where $x$ can be any real so that $\omega_1^x = \sigma$. In general, $\SCRL^{\Phi(x)}_\alpha$ and $T^{\Phi(x)}_\alpha$ are elements of $L_\sigma[\Phi(x)]$ when $\alpha < \sigma$ and is $\Delta_1$ definable in $L_\sigma[\Phi(x)]$ when $\alpha = \sigma$. The Solovay product lemma will indicate that each set belongs to $L_\sigma$ when $\alpha < \sigma$ and is $\Delta_1$ in $L_\sigma$ when $\alpha = \sigma$. 

The following will give a formal definition of $L^\sigma_\alpha$ and $T^\sigma_\alpha$ inside any model of $\KP$ and their basic properties.

\Begin{definition}{generic theory of a classification}
Let $z \in \cantorspace$. Let $\SCRL$ be a recursive langauge. Let $\Phi$ be a $\lborel(z)$ classification of $\F$ by $\SCRL$-structures. Let $\sigma$ be a countable $z$-admissible ordinal. Let $\CA$ be a countable model of $\KP$ containing $z$ such that $\mathrm{ON} \cap \CA = \sigma$.

Now work in $(L[z])^\CA$: Let $\bbS$ denote Steel's tagged tree forcing defined in $(L[z])^{\CA}$. Next, by $\Sigma_1$-recursion in $(L[z])^{\CA}$ (which is a model of $\KP$), define the function taking an ordinal $\alpha$ of $(L[z])^{\CA}$ to $\SCRL^\sigma_\alpha$ and $T^\sigma_\alpha$ as follows:
\begin{enumerate}[$\bullet$]
\item Let $\SCRL^\sigma_0 = \SCRL_{\omega\omega}$.

\item If $\alpha$ is limit ordinal, then let $\SCRL^\sigma_\alpha = \bigcup_{\beta < \alpha} \SCRL^\sigma_\beta$. 

\item For any $\alpha$, let $T^\sigma_\alpha$ be that set such that $1_\bbS \forces_\bbS$ ``$\check T_\alpha^\sigma$ is the complete theory of $\Phi(\dot T)$ in the fragment $\check \SCRL_\alpha^\sigma$''. Here $\dot T$ refers to the canonical name for the generic tree (construed as a real) produced by $\bbS$ as in Definition \ref{steel forcing}. Such a set exists using (I) of the next lemma to show $T^{\Phi(x)}_\alpha \in (L[z])^{\CA}$ for any $x \in (\cantorspace)^V$ such that $\omega_1^x = \sigma$. (Alternatively, $T^\sigma_\alpha$ is also the set of sentences $\varphi \in \SCRL^\sigma_\alpha$ so that $1_\bbS \forces \Phi(\dot T) \models \check \varphi$.)

\item Let $\SCRL_{\alpha + 1}^\sigma$ be the smallest fragment $\mathcal{F}$ of $\SCRL_{\infty\omega}$ so that $1_\bbS \forces_{\bbS}$ ``$\check{\mathcal{F}}$ is the smallest fragment of $\SCRL_{\infty\omega}$ so that for all $n \in \omega$, $\bigwedge p$ belongs to $\mathcal{F}$ for any nonprincipal $p \in S_n^{\SCRL^\sigma_\alpha}(T^\sigma_\alpha)$ realized by some $n$-tuple in $\Phi(\dot T)$''.
\end{enumerate}

\end{definition}

\Begin{lemma}{definability of theory of classification}
Assume the setting of Definition \ref{generic theory of a classification}. 
\begin{enumerate}[(I)]
\item For any $x \in \cantorspace \cap \CA$ such that $\omega_1^x = \sigma$ (in the real world) and any $\alpha \in \mathrm{ON}^\CA$, $\SCRL^{\Phi(x)}_\alpha, T^{\Phi(x)}_\alpha \in (L[z])^\CA$. Any $p \in S_n^{\SCRL_\alpha^\Phi(x)}(T^\Phi(x)_\alpha)$ realized by some tuple in $\Phi(x)$ belongs to $(L[z])^\CA$. 

\item The functions $\alpha \mapsto \SCRL^\sigma_\alpha$ and $\alpha \mapsto T^\sigma_\alpha$ are $\Delta_1$ in $(L[z])^\CA$. 

\item Let $T_\infty^\sigma = \bigcup_{\alpha \in \mathrm{ON}} T^\sigma_\alpha$ and $\SCRL_\infty^\sigma = \bigcup_{\alpha \in \mathrm{ON}} \SCRL^\sigma_\alpha$. Both are $\Delta_1$ classes in $(L[z])^\CA$. 

\item For any $x \in (\cantorspace)^V$ with $\omega_1^x = \sigma$ and $\alpha \leq \sigma$, $L^{\Phi(x)}_\alpha = L^\sigma_\alpha$ and $T^{\Phi(x)}_\alpha = T^\sigma_\alpha$. 
\end{enumerate}
\end{lemma}

\begin{proof}
(I) is proved using Lemma \ref{solovay product lemma}. 

(II) and (III) are proved much like the corresponding facts for $\SCRL^{\Phi(x)}_\alpha$ and $T^{\Phi(x)}_\alpha$ mentioned in Definition \ref{theory of model and Scott rank} and using the definability of the forcing relation $\forces_{\bbS}$. 

(IV) is proved by induction.
\end{proof}

\Begin{fact}{isolating formula T is pi}
Assume the setting of Lemma \ref{generic theory of a classification}. The relation ``$\varphi$ is an isolating formula of $S_n^{\SCRL_\infty^\sigma}(T_\infty^\sigma)$'' with free variable $\varphi$ ranging over $\SCRL^\sigma_\infty$ is a $\Pi_1$ relation.
\end{fact}

\begin{proof}
Recall that $T^\sigma_\infty$ and $\SCRL^\sigma_\infty$ are $\Delta_1$ and that $T^\sigma_\infty = \bigcup_{\alpha \in \mathrm{ON}} T^\sigma_\alpha$, where each $T^\sigma_\alpha$ is a complete theory in $\SCRL^\sigma_\alpha$. 

That $\varphi$ is an isolating formula can be expressed by saying for all $\psi \in \SCRL^\sigma_\infty$, for all $\beta \in \mathrm{ON}$, if $\varphi,\psi \in \SCRL^\sigma_\beta$, then either $(\forall \bar{x})(\varphi \Rightarrow \psi) \in T^\sigma_\beta$ or $(\forall \bar{x})(\varphi \Rightarrow \neg\psi) \in T^\sigma_\beta$. This can be formalized as a $\Pi_1$ statement.
\end{proof}

\Begin{fact}{classified by high scott rank}
(\cite{The-Countable-Admissible-Ordinal-Equivalence-Relation-pub} Theorem 4.2) Let $z \in \cantorspace$. Let $\SCRL$ be a recursive language. Suppose $\Phi : \cantorspace \rightarrow S(\SCRL)$ is a $\Delta_1^1(z)$ classification of $\F$ be $\SCRL$-structures. Then for all $x$ such that $\omega_1^x$ is $z$-admissible, $\SR(\Phi(x)) \geq \omega_1^x$. 
\end{fact}

\Begin{theorem}{classification has model high but not highest scott rank}
Let $\SCRL$ be a recursive language. Let $\Phi : \cantorspace \rightarrow S(\SCRL)$ be a $\lborel(z)$ classification of $F_{\omega_1}$ by $\SCRL$-structures. Let $\sigma$ be an ordinal so that $L_\sigma[z] \models \Sigma_2\text{-}\KP$ and $L_\sigma[z] \prec_1 L_{\omega_1^{L[z]}}[z]$. Let $x \in \cantorspace$ be such that $\omega_1^x = \sigma$. Then $\SR(\Phi(x)) = \sigma + 1$. 
\end{theorem}

Such a countable ordinal $\sigma$ can be found as follows: By the L\"owenheim-Skolem theorem, let $M \prec_\omega L_{\omega_1^{L[z]}}[z]$ be a countable elementary substructure containing $z$. Since $L_{\omega_1^{L[z]}}[z]$ thinks the transitive closure of all sets are countable, $M$ has this property too. For all $x \in M$, there is some bijection $f : \omega \rightarrow \mathrm{tc}(x)^M$, where $\mathrm{tc}$ denotes the transitive closure. Since $M$ is elementary, this $f$ really is a bijection of $\omega$ with $\mathrm{tc}(x)$. For all $n \in \omega$, $f(n) \in M$ so $\mathrm{tc}(x) \subseteq M$. This shows that $M$ is transitive. Since $M$ is a countable transitive elementary substructure of $L_{\omega_1^{L[z]}}[z]$, there is some countable $\sigma$ so that $M = L_{\sigma}[z]$ by condensation. Finally, $L_\sigma[z] \models \mathsf{ZF - P}$ (and in particular $\Sigma_2$-$\KP$) because $L_{\omega_1^{L[z]}}[z] \models \mathsf{ZF - P}$. 

\begin{proof}
Before beginning the proof, an outline will be given: For simplicity throughout the proof, suppose $\Phi$ is $\lborel$. By Lemma \ref{definability of theory of classification}, $\langle \SCRL^\sigma_\alpha : \alpha \in \mathrm{ON}\rangle$ and $\langle T^\sigma_\alpha : \alpha \in \mathrm{ON}\rangle$ are $\Delta_1$-classes in the constructible universe of any model of $\KP$ whose collection of standard ordinals has ordertype $\sigma$. In particular, these two sequences are $\Delta_1$-definable in $L_\sigma$. First, one will find an illfounded model of $\KP$, $\CB$, so that $\CB \cap \mathrm{ON} = \sigma$ and $\CB$ contains some real $c$ so that $\omega_1^c = \sigma$. Lemma \ref{definability of theory of classification} asserts that each $T^{\Phi(c)}_\alpha = T^\sigma_\alpha$ and $\SCRL^{\Phi(c)}_\alpha = \SCRL^\sigma_\alpha$, and hence they belong to $L_\sigma$. Furthermore, $\CB$ has the crucial property that any isolating formula for $\SCRL^\sigma_\infty$ in $L_\sigma$ is an isolating formula for $(\SCRL^\sigma_\infty)^\CB$ in $\CB$. This fact will be accomplished by simply arranging that $L_\sigma \prec_1 \CB$. This $\CB$ is found using Jensen's model existence theorem with an appropriate theory in a countable admissible fragment of $L_\sigma$ that attempts to express $\Sigma_1$-elementarity. The choice of the ordinal $\sigma$ will show that the theory to which the Jensen's model existence theorem is applied is consistent since it will be modeled by $L_{\omega_1^L}$. The purpose of the two sequences, $\langle \SCRL_\alpha^\sigma : \alpha \in \mathrm{ON}\rangle$ and $\langle T_\alpha^\sigma : \alpha \in \mathrm{ON}\rangle$, and the effort to establish that they belong to or are definable in the constructible universe of the relevant models of $\KP$ is to be able to express the absoluteness of being an isolating formula within this admissible fragment of $L_\sigma$.

By Fact \ref{classified by high scott rank}, $\SR(\Phi(c))$ is $\sigma$ or $\sigma + 1$. Now suppose for a contradiction that $\SR(\Phi(c)) = \sigma$. Within $L_\sigma[\Phi(c)]$, define a function $\Psi$ which assigns each tuple of $\Phi(c)$ to the least ordinal $\alpha$ so that there is some formula $\varphi \in \SCRL_\alpha^\sigma$ which is isolating and is realized by this tuple. This function is well-defined by the assumption that $\SR(\Phi(c)) = \sigma$. Now let $\tilde \Psi$ be the function defined in the same way but within the illfounded model $\CB$. Using the fact that $L_\sigma \prec_1 \CB$, the preservation of isolating formulas implies $\Psi = \tilde \Psi$. However, $\SR(\Phi(c)) = \sigma$ implies that the image of $\tilde \Psi$ is cofinal within the standard ordinals of $\CB$. Then by an overspill argument, $\tilde \Psi$ must take on some illfounded ordinal. This contradicts $\Psi = \tilde \Psi$. The details of the proof are given below.
\\*
\\*\indent Let $\mathscr{J}$ be the language consisting of the following objects:
\begin{enumerate}[(i)]
\item A binary relation symbol $\dot \in$.

\item New constant symbol $\hat{a}$ for each element of $a \in L_\sigma$. 

\item Two new constant symbols $\dot c$ and $\dot d$. 
\end{enumerate}
$\mathscr{J}$ can be considered a $\Delta_1$ definable language in the countable admissible set $L_\sigma$. 

Let $H$ be the theory in the countable admissible fragment $(\mathscr{J}_{\infty\omega})^{L_\sigma}$ with the following sentences:
\begin{enumerate}[(I)]
\item All the axioms of $\mathsf{ZF-P}$.

\item For each $a \in L_\sigma$, $(\forall v)(v \dot\in \hat{a} \Leftrightarrow \bigvee_{z \in a} \hat{z} = v)$.

\item Add the sentence ``$\dot c \in 2^{\hat{\omega}}$''. For each ordinal $\beta < \sigma$, add ``$L_{\hat\beta}[\dot c] \not\models \KP$''. 

\item For each $\Pi_1$ formula $\varphi(x_0, ..., x_{k - 1})$ of $\{\dot \in\}$ and elements $a_0, ..., a_{k - 1} \in L_\sigma$ such that $L_\sigma \models \varphi(a_0, ..., a_{k - 1})$, add the sentence ``$\varphi(\hat{a}_0, ..., \hat{a}_{k - 1})$''.

\item Add the sentence ``$\dot d$ is an ordinal''. For each $\beta < \sigma$, add ``$\dot d > \hat\beta$''. 
\end{enumerate}
The $\Sigma_1$-satisfaction relation of $L_\sigma$ is a $\Sigma_1$ relation in $L_\sigma$. (See \cite{the-fine-structure-of-the-constructible-hierarchy} Corollary 1.13.) Therefore, the $\Pi_1$-satisfaction relation of $L_\sigma$ is $\Pi_1$ definable in $L_\sigma$. So $H$ is a $\Pi_1$ and hence a $\Sigma_2$-definable subset of $L_\sigma$. Note that (III) states that each $\beta < \sigma$ is a not a $\dot c$-admissible ordinal. (IV) states that $L_\sigma$ will be a $\Sigma_1$ elementary substructure of any model of $H$. (V) states that $\dot d$ is an ordinal larger than each $\beta < \sigma$. 

(It will be seen below that one will only use the fact that isolating sentences in $L_\sigma$ remain isolating sentences in models of $H$. One can rewrite (IV) to express this rather than attempt to obtain full $\Sigma_1$-elementarity. Also (I), (II), (III) and an argument similar to the one below constitute Jensen's proof (see \cite{Admissible-Sets}) of Sack's theorem.)

$H$ is consistent: Let $\CB$ be the $\mathscr{J}$-structure with underlying domain $L_{\omega_1}$. Let $\dot \in^\CB = \in \upharpoonright L_{\omega_1}$. Let $c$ be any real in $L_{\omega_1^{L}}$ so that $\omega_1^c > \sigma$. Let $\dot c^\CB = c$. Let $\dot d = \sigma + 1$. For each $a \in L_{\sigma}$, let $\hat{a}^\CB = a$. (I), (II), (III), and (V) are clearly satisfied in $\CB$. Note that (IV) is satisfied since $L_\sigma \prec_1 L_{\omega_1^{L}}$. 

By the $\Sigma_2$ version of Jensen's model existence theorem (Fact \ref{jensen model existence theorem}), there is a $\SCRL$-structure $\CB$ so that $\CB \models H$, $\CB$ end extends $L_\sigma$, $\mathrm{WF}(\CB)$ is transitive, and $\mathrm{ON} \cap \CB = \sigma$. Let $c = \dot c^\CB$ and $d = \dot d^\CB$. By (V), $d$ is a nonstandard ordinal. So $\CB$ is an illfounded model. Note that $c \in \mathrm{WF}(\CB)$ since $c$ is a real.

First, to show that $\omega_1^c = \sigma$: For each $\beta < \sigma$, $\CB \models L_\beta[c] \not \models \KP$. Satisfaction is $\Delta_1$ so by absoluteness, $\mathrm{WF}(\CB) \models L_\beta[c] \not\models \KP$. Again by absoluteness, $V \models L_\beta[c] \not\models \KP$. Hence $\beta$ is not a $c$-admissible ordinal. This shows that $\omega_1^c \geq \sigma$. Also by (I), $\CB \models \mathsf{ZF - P}$. In particular, $\CB \models \KP$. By the truncation lemma (Fact \ref{truncation lemma}), $\mathrm{WF}(\CB) \models \KP$. Hence $\mathrm{WF}(\CB)$ is an admissible set containing $c$. So $\mathrm{ON} \cap \mathrm{WF}(\CB) = \sigma$ is a $c$-admissible ordinal. This shows $\sigma$ is the smallest $c$-admissible ordinal. By definition, $\omega_1^c = \sigma$. 

Since $\Phi$ is $\Delta_1^1$, $\Phi(c)$ is $\Delta_1^1(c)$. $\Phi(c)$ belong to any admissible set containing $c$. Thus $\Phi(c) \in \mathrm{WF}(\CB) \subseteq \CB$. 

By Fact \ref{classified by high scott rank}, $\SR(\Phi(c)) \geq \sigma$. Suppose toward a contradiction that $\SR(\Phi(c)) = \sigma$. Then $\Phi(c)$ is an atomic model of $T^{\Phi(c)}_{\sigma} = (T_\infty^\sigma)^{L_\sigma}$.

Since $\Phi(c) \in S(\SCRL)$, $\Phi(c)$ is an $\SCRL$-structure with underlying domain $\omega$. In $L_\sigma[\Phi(c)]$, define the function $\Psi : {}^{<\omega}\omega \rightarrow \mathrm{ON}$ by letting $\Psi(\bar{a})$ be the least ordinal $\alpha$ so that there is some $\varphi \in \SCRL^\sigma_\alpha$ with $\varphi$ an isolating formula for $T_\infty^\sigma$ in the fragment $\SCRL^\sigma_\infty$ realized by the tuple $\bar{a}$. 

Note that this is a well-defined function since $\Phi(c)$ is an atomic model of $(T_\infty^\sigma)^{L_\sigma[\Phi(c)]}$. (Also note that since $L_\sigma$ and $L_\sigma[\Phi(c)]$ have the same ordinals and hence the same constructible universe, $T^\sigma$ and $\SCRL^\sigma$ are the same class whether relativized in $L_\sigma$ or relativized in $L_\sigma[\Phi(c)]$.)

Now let $\tilde \Psi$ be the function defined in $\CB$ by the same formula used to defined $\Psi$ in $L_\sigma[\Phi(c)]$. 

(Note that it is not immediately seen that $\tilde \Psi$ is the same function as $\Psi$ since $\mathrm{ON}^\CB \neq \sigma = \mathrm{ON}^{L_\sigma[\Phi(c)]}$. Hence $(T^\sigma_\infty)^{\CB} \neq (T^\sigma_\infty)^{L_\sigma[\Phi(c)]}$ and $(\SCRL_\infty^\sigma)^\CB \neq (\SCRL_\infty^\sigma)^{L_\sigma[\Phi(c)]}$. In particular, it is not immediate that a formula in $(\SCRL_\infty^\sigma)^{L_{\sigma}[\Phi(c)]} = (\SCRL_\infty^\sigma)^{L_\sigma}$ which isolates a type in $L_\sigma$ would still isolate a type in the larger fragment $(\SCRL_\infty^\sigma)^{\CB}$. $\Sigma_1$-elementarity will be used to resolve this.)

The claim is that $\Psi = \tilde \Psi$: To see this, let $\bar{a}$ be a finite tuple of natural number understood to be a tuple from $\Phi(c)$. Since $\Phi(c)$ is an atomic model of $(T_\infty^\sigma)^{L_\sigma[\Phi(c)]} = (T_\infty^\sigma)^{L_\sigma}$, there is some $\varphi \in (\SCRL_\infty^\sigma)^{L_\sigma}$ so that $[\varphi]^{T^\sigma_\infty}_{(\SCRL^\sigma)^{L_\sigma}} = \{\mathrm{tp}^{\Phi(c)}_{(\SCRL_\infty^\sigma)^{L_\sigma}}(\bar{a})\}$. $L_\sigma \models$ ``$\varphi$ is a isolating formula for $S_n^{\SCRL_\infty^\sigma}(T_\infty^\sigma)$''. By Fact \ref{isolating formula T is pi}, this statement is $\Pi_1$. Since $L_\sigma \prec_1 \CB$. $\CB \models$ ``$\varphi$ is an isolating formula of $S_n^{\SCRL_\infty^\sigma}(T_\infty^\sigma)$''. 

This shows that in $\CB$, $[\varphi]_{\SCRL_\infty^\sigma}^{T_\infty^\sigma} = \{\mathrm{tp}^{\Phi(c)}_{\SCRL^\sigma}(\bar{a})\}$. Hence $\CB \models \tilde\Psi(\bar{a}) \leq \Psi(\bar{a})$. Suppose $\tilde\Psi(\bar{a}) < \Psi(\bar{a})$. There is some $\alpha < \Psi(\bar{a})$ and some formula $\psi \in (\SCRL^\sigma_\alpha)^{\CB}$ so that $\CB\models$ ``$\psi$ is an isolating formula which isolates $\mathrm{tp}^{\Phi(c)}_{\SCRL^\sigma}(\bar{a})$". Since $\Psi(\bar{a}) < \sigma$, $(\SCRL_\alpha^\sigma)^{\CB} = (\SCRL_\alpha^\sigma)^{L_\sigma}$. Hence $\psi \in L_\sigma$. Then by downward absoluteness of $\Pi_1$ statements from $\CB$ to $L_\sigma[\Phi(c)]$, $L_\sigma[\Phi(c)] \models$ ``$\psi$ is an isolating formula for $\mathrm{tp}^{\Phi(c)}_{\SCRL^\sigma}(\bar{a})$''. This contradicts the definition of $\Psi$ in $L_\sigma[\Phi(c)]$. 

This shows that $\Psi = \tilde \Psi$. Now suppose that $\Psi[\finNaturalSequence] = \sigma$. Since $\CB \models \mathsf{ZF - P}$ (in particular the full replacement axiom), $\tilde \Psi[\finNaturalSequence]$ is a set in $\CB$. Therefore, $\sup \tilde\Phi[\finNaturalSequence] \in \mathrm{ON}^\CB$. Since $\tilde \Psi = \Psi$ and $\Psi[\finNaturalSequence] = \sigma$, one must have that $\sup \tilde\Psi[\finNaturalSequence]$ is a nonstandard ordinal greater than each $\beta < \sigma$. Hence there is some $b < \sup\tilde\Psi[\finNaturalSequence]$ so that $\CB \models \beta < b$ for all standard $\beta < \sigma$ with $b \in \tilde\Psi[\finNaturalSequence]$. Thus there is some $\bar{a}$ so that $\Psi(\bar{a}) < b = \tilde\Psi(\bar{a})$. This contradicts $\Psi = \tilde \Psi$. 

This shows that $\Psi[\finNaturalSequence] < \sigma$ which implies $\SR(\Phi(c)) < \sigma$. However, it was already noted that $\SR(\Phi(c)) \geq \sigma$. Contradiction. This shows that $\SR(\Phi(c)) = \sigma + 1$. Since $\Phi$ is a classification, for any $x$ with $\omega_1^x = \omega_1^c = \sigma$, $\Phi(x) \approx_\SCRL \Phi(c)$. Hence $\SR(\Phi(x)) = \sigma + 1$. This completes the proof of the theorem.
\end{proof}

The following is still open.

\Begin{question}{high but not highest anywhere}
If $\SCRL$ is a recursive language and $\Phi : \cantorspace \rightarrow S(\SCRL)$ is a $\lborel$ classification of $\F$ by $\SCRL$-structures, then for all $x \in \cantorspace$, is $\SR(\Phi(x)) = \omega_1^x + 1$?
\end{question}

Note that the lightface $\lborel$ is important in the phrasing of the question. The $\F$-class $\{x : \omega_1^x = \omega_1^\emptyset\}$ is $\lborel(z)$ for any $z$ such that $\omega_1^z > \omega_1^\emptyset$. Therefore, with access to such a parameter $z$, one can easily modify a known $\lborel$ classification of $\F$ to obtain a $\lborel(z)$ classification that sends all elements of $\{x : \omega_1^x = \omega_1^\emptyset\}$ to some fixed recursive structure of Scott rank $\omega_1^\emptyset$ and leaves the other classes alone.

It seems that if a classification $\Phi$ is $\Delta_1^1(z)$, then the relativization of the above question should be to ask the same question but only for those $x \in \cantorspace$ such that $\omega_1^x$ is $z$-admissible.

\bibliographystyle{amsplain}
\bibliography{references}

\end{document}